\documentclass[12pt,a4paper]{article}
\usepackage[T2A]{fontenc}
\usepackage[english]{babel}
\usepackage{latexsym,amsfonts,amssymb,amsmath,longtable,mathrsfs,amsthm}
\usepackage[unicode,final,hyperindex]{hyperref}
\renewcommand{\le}{\leqslant}
\renewcommand{\ge}{\geqslant}

 \DeclareMathOperator{\Aut}{Aut}
\DeclareMathOperator{\Sym}{Sym}

\DeclareMathOperator{\Alt}{Alt}

\def\@begintheorem#1#2{\trivlist \item[\hskip \labelsep%
{\hspace*{\parindent}\bf #1\ #2{.}}]\it}
\def\@opargbegintheorem#1#2#3{\trivlist
      \item[\hskip \labelsep{\hspace*{\parindent}\bf #1\ #2\ \rm #3{.}}]\it}

\newtheorem{teo}{Theorem}
\newtheorem{lem}{Lemma}
\newtheorem{hyp}{Conjecture}

\newtheorem{cor}{Corollary}
\renewcommand{\proof}{\par\mbox{P r o o f.}\ \ }

\newenvironment{Titul}{\begin{center}}{\end{center}
\baselineskip 12pt
\addtolength{\textheight}{-14mm}
}

\newenvironment{Anot}{\begin{quote}\scriptsize\mbox{~~~}}{\end{quote}}

\begin{document}

\begin{Titul}
On the intersections of solvable Hall subgroups in finite groups\footnote{
The work is supported by SB RAS and UrB RAS, Integration project ``Groups and graphs''. The first author is supported by RF President grant for young scientists, MK-3036.2007.1, RF President grant for scientific schools NSc-344.2008.1, and RFBR, grant N  08-01-00322. The second author is supported by RFBR, grant N 07-01-00148.}
\end{Titul}

\begin{Anot}
In the paper we consider the following conjecture: if a finite group  $G$ possesses a solvable $\pi$-Hall subgroup $H$, then there exist elements  $x,y,z,t\in G$ such that the identity $H\cap H^x\cap H^y\cap H^z\cap H^t=O_\pi(G)$ holds. The minimal counter example is shown to be an almost simple group of Lie type.

\end{Anot}

\section*{Introduction}

The notation in the paper is standard and agree with that of  \cite{ATLAS}. In this paper we consider finite groups only, so the term  ``group'' always means ``finite group''. By $\pi$ we always denote a set of primes, by $\pi'$ its complement in the set of all primes is denoted. If   $n$ is a positive integer, then a set of all prime divisors of $n$ is denoted by $\pi(n)$. A number $n$ is called a {\em $\pi$-number}, if $\pi(n)\subseteq \pi$. For a finite group $G$ we set  $\pi(G):=\pi(\vert G\vert)$. A group $G$ is called a {\em $\pi$-group}, if $\pi(G)\subseteq \pi$. The maximal normal  $\pi$-subgroup of  $G$ is denoted by $O_\pi(G)$. A subgroup $H$ of $G$ is called a {\em $\pi$-Hall} subgroup, if $\pi(H)\subseteq\pi$ and $\pi(\vert
G:H\vert)\subseteq \pi'$. The set of all $\pi$-Hall subgroups of $G$ is denoted by $\mathrm{Hall}_\pi(G)$, the set of all Sylow $p$-subgroups of $G$ is denoted by $\mathrm{Syl}_p(G)$ (clearly, $\mathrm{Syl}_p(G)=\mathrm{Hall}_p(G)$). A group $G$ is called a {\em $\pi$-solvable}, if it possesses a subnormal series $\{e\}=G_0<G_1<G_2<\ldots<G_{n-1}<G_n=G$ such that all its section are either $\pi'$-groups, or solvable groups. For a group $G$ and a subgroup $S$ of $S_n$ by $G\wr S$ a permutation wreath product is denoted.

After Ph.Hall \cite{Hall} we say that a group  $G$ {\it satisfies $E_\pi$} (or briefly $G\in E_\pi$), if $G$ possesses a  $\pi$-Hall subgroup. If
$G\in E_\pi$ and all $\pi$-Hall subgroups are conjugate, then we say that $G$ {\it satisfies $C_\pi$} ($G\in
C_\pi$). If $G\in C_\pi$ and every $\pi$-subgroup of $G$ is included into a $\pi$-Hall subgroup of $G$, then we say that $G$ {\it satisfies $D_\pi$} ($G\in
D_\pi$).

Let $A,B,H$ be subgroups of $G$ such that $B\unlhd A$. Then $N_H(A/B)=N_H(A)\cap N_H(B)$ is called a normalizer of $A/B$. If $x\in N_H(A/B)$, then $x$ induces an automorphism of  $A/B$, acting by $Ba\mapsto B x^{-1}ax$. Thus there exists a homomorphism $N_H(A/B)\rightarrow \text{Aut}(A/B)$. The image of $N_H(A/B)$ under this homomorphism is denoted by $\text{Aut}_H(A/B)$ and is called a 
{\em group of induces automorphisms} of~${A/B}$ in $H$. Note that for given composition factor $A/B$ of $G$ the group of induced automorphisms  $\Aut_G(A/B)$ depends on the choice of $A$ and $B$, i.~e., depends on the choice of a composition series. If  $A\leq G$, then $\Aut_G(A)=\Aut_G(A/\{e\})$ by definition.

In this paper we consider the following 

\begin{hyp}\label{Mainhypecture}
Let $H$ be a solvable  $\pi$-Hall subgroup of a finite group $G$. Then there exist elements $x,y,z,t$ such that the equality
\begin{equation}\label{MainIdentity}
H\cap H^x\cap H^y\cap H^z\cap H^t=O_\pi(G)
\end{equation}
holds
\end{hyp}

We show that a minimal counter example to this conjecture is an almost simple group, and also we show that Conjecture \ref{Mainhypecture} is satisfied for almost simple groups with simple socle isomorphic to an alternating or a sporadic group.

Let $\overline{\phantom{G}}:G\rightarrow G/O_\pi(G)$ be the natural homomorphism. Using elementary equation $\vert A\cdot
B\vert=\frac{\vert A\vert\cdot\vert B\vert}{\vert A\cap B\vert}$, where $A$ and $B$ are subgroups of  $G$, we obtain (under Conjecture \ref{Mainhypecture}) that   for a $\pi$-Hall subgroup $H$ the series of inequalities

\begin{multline*}
\vert \overline{G}\vert\ge \vert\overline{H}\cdot \overline{H}^x\vert=\frac{\vert \overline{H} \vert\cdot\vert
\overline{H}^x\vert}{\vert \overline{H}\cap \overline{H}^x\vert}= \frac{\vert \overline{H} \vert\cdot\vert
\overline{H}^x\vert}{\vert \overline{H}\cap \overline{H}^x\vert}\cdot \frac{\vert \overline{H}\cap \overline{H}^x\vert\cdot\vert
\overline{H}^y\vert}{\vert \overline{H}\cap \overline{H}^x\cap \overline{H}^y\vert\cdot \vert (\overline{H}\cap
\overline{H}^x)\overline{H}^y\vert} \ge\\ \frac{\vert \overline{H}\vert\cdot\vert \overline{H}^x\vert\cdot\vert
\overline{H}^y\vert}{\vert \overline{H}\cap \overline{H}^x\cap \overline{H}^y\vert\cdot\vert \overline{G}\vert}\ge\ldots\ge
\frac{\vert \overline{H} \vert^5}{\vert  \overline{G}\vert^{3}}
\end{multline*}
holds. Thus an affirmative answer to Conjecture \ref{Mainhypecture} implies that the inequality $\vert H/O_\pi(G)\vert<\vert
G:H\vert^4$ holds.

Note that D.S.Passman in \cite{Passman} proved that a $p$-solvable group always possesses three Sylow  $p$-subgroups such that their intersection is equal to  $O_p(G)$. Later V.I.Zenkov proved the same statement for an arbitrary group (see \cite[Corollary~C]{ZenkovSylow}). In \cite{DolfiOdd} S.Dolfi proved that if $2\not\in\pi$, then every $\pi$-solvable group  $G$ possesses three $\pi$-Hall subgroups such that their intersection is equal to  $O_\pi(G)$. In \cite{Dolfi} S.Dolfi proved that every in $\pi$-solvable group $G$ there exist elements $x,y\in G$ such that
the equality  $H\cap H^x\cap H^y=O_\pi(G)$ holds (see also~\cite{VdovinIntersSolv}). On the other hand, the condition of solvability of  $H$ in Conjecture \ref{Mainhypecture} is essential. Indeed, if $p$ is a prime, then  $S_{p-1}$ is a $p'$-Hall subgroup of  $S_p$, but the intersection of every   $p-2$ conjugate subgroups is not equal to~$\{e\}$.

Let  $G$ be a subgroup of the symmetric group $\Sym_n$. A partition $P_1\sqcup P_2\sqcup\ldots\sqcup P_m$ of $\{1,\ldots,n\}$ is called an 
{\em asymmetric partition} for  $G$, if only the identity element of $G$ fixes the partition, i.~e., the equality $P_jx=P_j$ for all
$j=1,\ldots,m$ implies that $x=e$. Clearly for every subgroup $G$  the partition $P_1=\{1\},
P_2=\{2\},\ldots,P_n=\{n\}$ is always asymmetric. In  \cite[Theorem~1.2]{Seress} is proven that for $G$ solvable there exists an asymmetric partition with~${m\le5}$.

\section{Preliminary results}

The following statements are known.

\begin{lem}\label{base}
Let $A$ be a normal subgroup of $G$. Then the following statements hold:
\begin{itemize}
\item[{\em (a)}] for $H$ a $\pi$-Hall subgroups of $G$ we have that $HA/A$ and $H\cap A$ are $\pi$-Hall subgroups of $G/A$ and
$A$ respectively;
\item[{\em (b)}] if $A\in C_\pi$ and $G/A\in E_\pi$ (resp. $G/A\in C_\pi$), then $G\in E_\pi$ (resp. $G\in C_\pi$);
\item[{\em (c)}] if there exists a subnormal series of $G$ such that all its factors are either  $\pi$- or $\pi'$- groups, then 
$G\in D_\pi$.
\item[{\em (d)}]  if $G/A$ is a $\pi$-group and $H\in \mathrm{Hall}_\pi(A)$, then a $\pi$-Hall subgroup $\overline{H}$ of $G$ with
 $\overline{H}\cap A=H$ exists if and only if  $G$ acting by conjugation leaves the set $\{H^a\mid a\in A\}$ invariant.
\end{itemize}
\end{lem}

Let $S$ be a non-Abelian finite simple group and   $G$ is such that there exists a normal subgroup $T=S_1\times\ldots\times
S_n$ of $G$ satisfying the following conditions
\begin{itemize}
\item[(a)] $S_1\simeq\ldots\simeq S_n\simeq S$;
\item[(b)] the action of $G$ on $\{S_1,\ldots,S_n\}$ by conjugation is transitive;
\item[(c)] $C_G(S_1,\ldots,S_n)=\{e\}$.
\end{itemize}

Subgroups $N_G(S_1),\ldots,N_G(S_n)$ are conjugate, since the action of $G$ is transitive. Let $\rho:G\rightarrow \Sym_n$ be the permutation representation of $G$ on the right cosets of $N_G(S_1)$. Since the action by right multiplication of $G$ on the right cosets of  $N_G(S_1)$ coinside with the action by conjugation of $G$ on the set $\{S_1,\ldots,S_n\}$ we obtain that $G\rho$ is a transitive subgroup of  $\Sym_n$. By \cite[Hauptsatz~1.4, p.~413]{Huppert}
there exists a monomorphism $$\varphi:G\rightarrow (N_G(S_1)\times \ldots\times N_G(S_n)): (G\rho)=N_G(S_1)\wr (G\rho)=L.$$
Consider the natural homomorphism  $$\psi:L\rightarrow L/(C_G(S_1)\times\ldots\times C_G(S_n)).$$ Denoting the subgroup
$\Aut_G(S_i)=N_G(S_i)/C_G(S_i)$ by $A_i$ we obtain that  $$\varphi\circ\psi:G\rightarrow (A_1\times\ldots\times
A_n):(G\rho)$$ is a homomorphic inclusion of $G$ into  $(A_1\times\ldots\times A_n):(G\rho)\simeq A_1\wr (G\rho)=:\overline{G}$. The kernel of the homomorphism is equal to  $C_G(S_1,\ldots,S_n)=\{e\}$, i.~e., $\varphi\circ\psi$ is a monomorphism and we identify  $G$ with the subgroup $G(\varphi\circ\psi)$ of~$\overline{G}$.

\begin{lem}\label{NonAbelianSimpleDirectAutomorphisms}
In the above notation assume that there exists a  $\pi$-Hall subgroup $H$ of $G$ and $G=(S_1\times\ldots\times S_n)H$. Then there exists a  $\pi$-Hall subgroup  $\overline{H}$ of $\overline{G}$ such that $\overline{H}\cap G=H$.
\end{lem}

\proof
The factor group $\overline{G}/(S_1\times\ldots\times S_n)$ is a  $\pi$-group. By construction $\Aut_{\overline{G}}(S_1)=\Aut_G(S_1)$ and $\Aut_{\overline{G}}(S_i)=\Aut_G(S_i)\simeq \Aut_G(S_1)$ for all $i$.  By Lemma \ref{base}(d) a $\pi$-Hall subgroup  $M$ of $S_1\times \ldots\times S_k$ is included into a  $\pi$-Hall subgroup of  $G$ (resp. of $\overline{G}$) if and only if the set $\{M^s\mid s\in S_1\times\ldots\times S_k\}$ is invariant under the action by conjugation of  $G$ (resp. of $\overline{G}$). In order to complete the proof it is enough to show that for $H$ a $\pi$-Hall subgroup  of   $G$ and  $M=H\cap(S_1\times\ldots\times S_n)$ the class $\{M^s\mid s\in S_1\times\ldots\times
S_n\}$ is $\overline{G}$-invariant. Indeed, in this case there exists a   $\pi$-Hall subgroup   $\overline{H}$ of $\overline{G}$ with $\overline{H}\cap (S_1\times\ldots\times S_n)=M$. Hence $H,\overline{H}\leq N_{\overline{G}}(M)$. Since $\overline{G}/(S_1\times\ldots\times S_n)$ is a $\pi$-group Lemma  \ref{base}(c) implies that $N_{\overline{G}}(M)\in D_\pi$, so $H$ is conjugate to a subgroup of $\overline{H}$. The fact that  $\{M^s\mid s\in S_1\times\ldots\times S_n\}$ is $\overline{G}$-invariant follows easily from the construction of $\overline{G}$ and inclusion of $G$ into $\overline{G}$, and also from the fact that  $\{M^s\mid s\in S_1\times\ldots\times S_n\}$ is $G$-invariant by Lemma \ref{base}(d). Indeed, every element from $\overline{G}$ can be written as $(a_1,\ldots,a_n)g$, where $a_i\in A_i$ and $g\in G$. Since  $\{M^s\mid s\in S_1\times\ldots\times S_n\}$ is
$G$-invariant, we may assume that  $g=e$. More over, for every $i$ the class  $\{(M\cap S_i)^x\mid x\in S_i\}$ is
$A_i$-invariant by construction, an the lemma follows.
\qed

\begin{lem} \label{IntersectionAbelianSubgroups} {\em \cite{ZenkovAbelian}}
Let $A$ be an Abelian subgroup of a finite group $G$. Then there exists $x\in G$ such that $A\cap A^x\leq F(G)$.
\end{lem}

\begin{cor}\label{IntersectionAbelHallsubgroups}
Let $G=P: R$, where $P$ is a normal Sylow $p$-subgroup of  $G$ and $R$ is an Abelian  $p'$-Hall subgroup of $G$ with $C_R(P)=\{e\}$. Then there exists $x\in P$ such that~${R\cap R^x=\{e\}}$.
\end{cor}

\section{Reduction to an almost simple group}

For inductive arguments we need to define an additional condition {\bfseries(Orb)} to Conjecture \ref{Mainhypecture}, which also need a verification. It is clear that for $h\in H$ and four-tuple $(H^x,H^y,H^z,H^t)$ satisfying \eqref{MainIdentity}, the four-tuple 
$(H^{xh},H^{yh},H^{zh},H^{th})$ also satisfies \eqref{MainIdentity}. Thus  $H$ acts by conjugation on the set of four-tuples  $(H^x,H^y,H^z,H^t)$ satisfying \eqref{MainIdentity}. Condition {\bfseries(Orb)} states that  there exists at least five orbits under this action of~$H$.

Let $G$ be a finite group and subgroups $A,B$ of $G$ are chosen so that  $A/B$ is a non-Abelian composition factor of   $G$. The group $G$ is said to satisfy condition $\text{{\bfseries(CI)}}_\pi$, if for every composition factor  $A/B$ (and for every choice of $A,B$) Conjecture \ref{Mainhypecture} and condition {\bfseries(Orb)} are satisfied for~${\Aut_G(A/B)}$.

\begin{teo}\label{Mainteorem}
Let $H$ be a solvable $\pi$-Hall subgroup of $G$, and $G$ satisfies  {\em $\text{{\bfseries(CI)}}_\pi$}.  Then there exist $x,y,z,t\in G$ such that {\em\eqref{MainIdentity}} holds. Moreover, if $G$ is insolvable, then it satisfies~{\em{\bfseries(Orb)}}.
\end{teo}

\proof
Let $G$ be a counter example to the statement of the theorem of minimal order. By \cite{Dolfi} (see also \cite{VdovinIntersSolv})  $G$ is insolvable. Let $S(G)$ be a solvable radical of $G$. Note that solvability of $H$ implies that $O_\pi(G)$ is also solvable, so~$O_\pi(G)\leq S(G)$.

Assume that $S(G)$ is nontrivial. By Lemma \ref{base}(a),  $H_1=H\cap S(G)$ is a $\pi$-Hall subgroup of $S(G)$ and
$\vert S(G)\vert<\vert G\vert$. Hence $S(G)$ satisfies to the statement of the theorem. More over, denoting $G/S(G)$ by $\overline{G}$ we obtain that  $\overline{G}$ also satisfies to the statement of the theorem, and
$O_\pi(\overline{G})=\{\overline{e}\}$. The Frattini argument  (since $S(G)$ is solvable, all its $\pi$-Hall subgroups are conjugate) implies that the equality $G=N_G(H_1)S(G)$ holds. Choose elements    $x_1,y_1,z_1,t_1\in S(G)$ so that the equality 
$$H_1\cap H_1^{x_1}\cap H_1^{y_1}\cap H_1^{z_1}\cap
H_1^{t_1}=O_\pi(S(G))=O_\pi(G)$$ is true. Elements  $x_2,y_2,z_2,t_2\in G$ we choose so that the equality 
\begin{equation}\label{first}
\overline{H}\cap
\overline{H}^{\bar{x}_2}\cap \overline{H}^{\bar{y}_2} \cap \overline{H}^{\bar{z}_2} \cap \overline{H}^{\bar{t}_2}=\{\bar{e}\}
\end{equation} is true. Now choose a $\pi$-Hall subgroup $H$ in $N_G(H_1)$, in $N_G(H_1^{x_1})$ we choose a $\pi$-Hall subgroup $H^x$ so that its image in 
$\overline{G}$ coincides with  $\overline{H}^{\bar{x}_2}$ and so on; in  $N_G(H_1^{t_1})$ we choose a $\pi$-Hall subgroup  $H^t$ so that its image in $\overline{G}$ coincides with  $\overline{H}^{\bar{t}_2}$ (such subgroups exist by Lemma~\ref{base}(d)).

We show that  $D=H\cap H^x\cap H^y\cap H^z\cap H^t=O_\pi(G)$. We have
\begin{multline*}
D\cap S(G)=H\cap H^x\cap H^y\cap H^z\cap H^t\cap S(G)=\\ (H\cap S(G))\cap (H^x\cap S(G))\cap (H^y\cap S(G))\cap (H^z\cap S(G))\cap (H^t\cap
S(G))=\\ H_1\cap H_1^{x_1}\cap H_1^{y_1}\cap H_1^{z_1}\cap H_1^{t_1}=O_\pi(G),
\end{multline*}
i.~e. $D\cap S(G)=O_\pi(G)$. More over
\begin{equation*}
\overline{D}= \overline{H}\cap
\overline{H}^{\bar{x}_2}\cap \overline{H}^{\bar{y}_2} \cap \overline{H}^{\bar{z}_2} \cap \overline{H}^{\bar{t}_2}=\{\bar{e}\},
\end{equation*}
so $D=O_\pi(G)$. Since   $G$ is insolvable, the factor group  $\overline{G}$ is insolvable as well. By induction, there exist at least five four-tuples $(\overline{H}^{\bar{x}_2},\overline{H}^{\bar{y}_2},\overline{H}^{\bar{z}_2},\overline{H}^{\bar{t}_2})$, satisfying \eqref{first}, under action by conjugation of $H$. Therefore there exists at least five four-tuples $(H^x,H^y,H^z,H^t)$ satisfying~\eqref{MainIdentity}. A contradiction with the fact, that $G$ is a minimal counter example. Hence~${S(G)=\{e\}}$.

Consider $G_1=F^\ast(G)H$. Then $C_G(F^\ast(G))\leq F^\ast(G)$, and, since $S(G)=\{e\}$, we have that $F(G)=\{e\}$ and
$F^\ast(G)=E(G)=S_1\times\ldots\times S_k$, where  $S_1,\ldots,S_k$ are non-Abelian simple groups. Therefore $H$ acts by conjugation faithfully on
 $E(G)$. If $G\not=G_1$, then $G_1$ satisfies to the statement of the theorem by induction. Hence $G$ also satisfies to the statements of the theorem, i.~e.,  $G$ is not a minimal counter example. So  $G=G_1=E(G)H$. Since $S_1,\ldots,S_k$ are simple we obtain that $G$ acting by conjugation permutes elements of $\{S_1,\ldots,S_k\}$.

Assume that $G$ or, equivalently  $H$ acts intransitively on the set $\{S_1,\ldots,S_k\}$. Since  $E(G)=S_1\times\ldots\times S_k$ we have that  $E_1=\langle S_1^H\rangle\not=E(G)$ is a normal subgroup of $G$. So   $E(G)=E_1\times E_2$, where $E_1$ and $E_2$ are $H$-invariant subgroups.  Therefore there exists a homomorphism   $G\rightarrow G/(C_G(E_1))\times G/(C_G(E_2))$ such that it image is a subdirect product of $G/(C_G(E_1))$ and $G/(C_G(E_2))$, while the kernel is  $C_G(E_1)\cap C_G(E_2)=C_G(E(G))=\{e\}$. Denote the projections of $G$ onto $G/(C_G(E_1))$ and $G/(C_G(E_2))$ by $\pi_1$ and $\pi_2$ respectively.  Since $G=E(G)H$ and $E_1\leq \mathrm{Ker}(\pi_2)$, $E_2\leq \mathrm{Ker}(\pi_1)$ we obtain the following equalities
$G\pi_1=E_1 (H\pi_1)$ and $G\pi_2=E_2 (H\pi_2)$ (we identify $E_i\pi_i$ and $E_i$, since  $E_i\pi_i\simeq E_i$). By induction there exist $x_i,y_i,z_i,t_i\in E_i(H\pi_i)$ such that
\begin{equation}\label{Proj}
(H\pi_i)\cap (H\pi_i)^{x_i}\cap (H\pi_i)^{y_i}\cap (H\pi_i)^{z_i}\cap
(H\pi_i)^{t_i}=\{e\}.
\end{equation}
Since $G\pi_i=E_i(H\pi_i)$ we may assume that elements $x_i,y_i,z_i,t_i$ are in  $E_i$. Consider
$x=x_1x_2$, $y=y_1y_2$, $z=z_1z_2$, $t=t_1t_2$. Since \eqref{Proj} is true for all $i$ it follows that for elements $x,y,z,t$ equality \eqref{MainIdentity} is true. Since there exist at least five orbits of four-tuples
$((H\pi_1)^{x_1},(H\pi_1)^{y_1},(H\pi_1)^{z_1},(H\pi_1)^{t_1})$ under the action of  $H\pi_1$ we obtain that there exist at least five four-tuples $(H^x,H^y,H^z,H^t)$ under the action of  $H$. Thus  $G$ satisfies to the statements of the theorem, i.~e.,  $G$ is not a counter example.

Thus $H$ acts transitively on $\{S_1,\ldots,S_k\}$ and, in view of condition $\text{{\bfseries(CI)}}_\pi$, we may assume that
$k>1$. By Lemma  \ref{NonAbelianSimpleDirectAutomorphisms} we also may assume that $G=(A_1\times\ldots\times A_k): L=A_1\wr L$, where
$S_i\leq A_i\leq \Aut(S_i)$ and $L$ is the image of  $G$ or, that is the same, of $H$ in $\Sym_k$ (in particular $L$ is a solvable 
$\pi$-group). Denote  $H\cap A_i$ by $H_i$. By Lemma \ref{base} the subgroup $H_i$ is a solvable  $\pi$-Hall subgroup of $A_i$. In view of condition $\text{{\bfseries(CI)}}_\pi$, there exist at least five orbits of four-tuples $(H_1^{x_1},H^{y_1},H^{z_1},H^{t_1})$, satisfying the equality
\begin{equation*}
H_1\cap H_1^{x_1}\cap H_1^{y_1}\cap H_1^{z_1}\cap H_1^{t_1}=\{e\},
\end{equation*} under the action of  $H_1$
Set $(H_1^{x_{1,j}},H_1^{y_{1,j}},H_1^{z_{1,j}},H_1^{t_{1,j}})$ to be a representative of the  $j$-th orbit ($j=1,2,3,4,5$). Let $h\in H$ and
$S_1^h=S_i$ (hence, $H_1^h=H_i$). Consider  $$(H_1^{x_{1,j}},H_1^{y_{1,j}},H_1^{z_{1,j}},H_1^{t_{1,j}})\text{ and }
(H_1^{x_{1,l}},H_1^{y_{1,l}},H_1^{z_{1,l}},H_1^{t_{1,l}})$$, representatives of distinct orbits of four-tuples under the action of  $H_1$. It is clear that $$((H_1^{x_{1,j}})^h,(H_1^{y_{1,j}})^h,(H_1^{z_{1,j}})^h,(H_1^{t_{1,j}})^h)\text{ and }
((H_1^{x_{1,l}})^h,(H_1^{y_{1,l}})^h,(H_1^{z_{1,l}})^h,(H_1^{t_{1,l}})^h)$$ are also representatives of distinct orbits of four-tuples under the action of $H_i$. More over, the set of four-tuples  $$\{((H_1^{x_{1,j}})^h\cap
S_i,(H_1^{y_{1,j}})^h\cap S_i,(H_1^{z_{1,j}})^h\cap S_i,(H_1^{t_{1,j}})^h\cap S_i)\mid h\in H\}$$ is in the same orbit under the action of
$H_i$. Note that the subgroup $H_{1,j_1}\times\ldots\times H_{k,j_k}$ is a $\pi$-Hall subgroup of  $S_1\times\ldots\times
S_k$ (see Lemma  \ref{base}(a)) and is contained in a  $\pi$-Hall subgroup of   $G$ for every choice of $j_i$-s from
$\{1,2,3,4,5\}$ (see Lemma~\ref{base}(d)).

By \cite[Theorem~1.2]{Seress} there exists an asymmetric partition $P_1\sqcup P_2\sqcup P_3\sqcup P_4\sqcup P_5$ множества $\{1,\ldots,k\}$
(some of these sets may be empty) such that only identity element of $L$ ($\simeq G/(S_1\times\ldots\times S_k)$) stabilizes this partition.
Set  $j(i)=m$, if $i\in P_m$ and consider subgroups 
\begin{multline*}
M_1=H_1\times\ldots\times H_k, M_2=H_1^{x_{1,j(1)}}\times \ldots\times H_k^{x_{k,j(k)}}, M_3=H_1^{y_{1,j(1)}}\times
\ldots\times H_k^{y_{k,j(k)}},\\ M_4=H_1^{z_{1,j(1)}}\times \ldots\times
H_k^{z_{k,j(k)}}, M_5=H_1^{t_{1,j(1)}}\times \ldots\times H_k^{t_{k,j(k)}}.
\end{multline*}
By construction  there exist elements  $x,y,z,t\in S_1\times\ldots\times S_k$ such that $M_2=M_1^x, M_3=M_1^y,M_4=M_1^z,M_5=M_1^t$.  We show that $H\cap H^x\cap H^y\cap H^z\cap H^t=\{e\}$. Note that by construction 
$H\cap (S_1\times\ldots\times S_k)=M_1$, $H^x\cap (S_1\times\ldots\times S_k)=M_2$, $H^y\cap (S_1\times\ldots\times S_k)=M_3$, $H^z\cap
(S_1\times\ldots\times S_k)=M_4$, $H^t\cap (S_1\times\ldots\times S_k)=M_5$, so if $h\in H\cap H^x\cap H^y\cap H^z\cap H^t$, then  $h$
normalizes subgroups  $M_1,M_2,M_3,M_4,M_5$. Assume that  $S_i^h=S_{h(i)}$, then $H_i^h=H_{h(i)}$, while tuples $$((H_i^{x_{i,j(i)}})^h, (H_i^{y_{i,j(i)}})^h,(H_i^{z_{i,j(i)}})^h,(H_i^{t_{i,j(i)}})^h), (H_{h(i)}^{x_{h(i),j(h(i))}},H_{h(i)}^{y_{h(i),j(h(i))}},H_{h(i)}^{z_{h(i),j(h(i))}},H_{h(i)}^{t_{h(i),j(h(i))}})$$
are in the same $H_{h(i)}$-orbit. Hence   $j(i)$ and $j(h(i))$ are in the same set  $P_j$ for $j=1,2,3,4,5$. Therefore  $h$ stabilizes the partition $P_1\sqcup P_2\sqcup P_3\sqcup P_4\sqcup P_5$ of  $\{1,\ldots,k\}$, so its image in $L$ equals $e$. Hence $h\in M_1\cap M_2\cap M_3\cap M_4\cap M_5=\{e\}$.

Chousing  $M_2,M_3,M_4,M_5$ to be equal to $H_1^{x_{1,j(1)}+1}\times \ldots\times H_k^{x_{k,j(k)}+1}$, $H_1^{y_{1,j(1)+1}}\times
\ldots\times H_k^{y_{k,j(k)+1}}$, $H_1^{z_{1,j(1)+1}}\times \ldots\times H_k^{z_{k,j(k)+1}}$, $H_1^{t_{1,j(1)+1}}\times \ldots\times
H_k^{t_{k,j(k)+1}}$,\ldots, $H_1^{x_{1,j(1)}+4}\times \ldots\times H_k^{x_{k,j(k)}+4}$, $H_1^{y_{1,j(1)+4}}\times
\ldots\times H_k^{y_{k,j(k)+4}}$, $H_1^{z_{1,j(1)+4}}\times \ldots\times H_k^{z_{k,j(k)+4}}$, $H_1^{t_{1,j(1)+4}}\times \ldots\times
H_k^{t_{k,j(k)+4}}$ we obtain at least five orbits of four-tuples, satisfying \eqref{MainIdentity}, under the action of~$H$.
\qed

\section{Intersection of solvable Hall subgroups in almost simple groups}

In this section we prove that a finite almost simple group, with simple socle isomorphic to either an alternating or a sporadic group, satisfies~$\text{{\bfseries(CI)}}_\pi$.

\begin{teo}\label{AlmostSimpleIntersections}
Let $H$ be a solvable  $\pi$-Hall subgroup of an almost simple group $G$. Assume also that the socle  $F^*(G)$ is isomorphic to either an alternating group, or a sporadic group. Then there exist  $x,y,z,t\in G$ such that equality {\em\eqref{MainIdentity}} holds. More over,  $G$ satisfies~{\em{\bfseries(Orb)}}.
\end{teo}

\proof
We show first that if $G$ possesses elements $x,y,z$ such that $H\cap H^x\cap H^y\cap H^z=\{e\}$, then there exists at least five orbits of four-tuples $(H^x,H^y,H^z,H^t)$, satisfying \eqref{MainIdentity}, under the action by conjugation of $H$. Indeed, if  $H\cap H^x\cap
H^y\cap H^z=\{e\}$, then four-tuples $$(H,H^x,H^y,H^z), (H^x,H,H^y,H^z), (H^x,H^y,H,H^z), (H^x,H^y,H^z,H), (H^x,H^x,H^y,H^z)$$
are in distinct $H$-orbits.

Now we proceed by investigating distinct almost simple groups. For given almost simple group we try to find elements  $x,y,z$, with $H\cap H^x\cap H^y\cap H^z=\{e\}$ first (as we noted above, in this case the condition on the number of $H$-orbits is satisfied automatically), and only if such elements do not exist, we find four elements $x,y,z,t$ and prove that there exist at least five $H$-orbits. Set $S=F^\ast(G)$ to be a simple socle of~$G$.

{\bfseries (I)} $S\simeq \mathrm{Alt}_n$, $n\ge5$. By \cite[Theorem~4.3]{RevVdoDpiFinal} and Lemma  \ref{base}(d) it follows that every $\pi$-Hall subgroup of  $S$ is contained in a $\pi$-Hall subgroup of $\Aut(S)$. Hence \cite[Theorem~А4]{Hall} implies that either $\vert\pi\cap\pi(S)\vert=1$ and $H$ is a Sylow subgroup of $S$ (hence of $G$), or $n=5$, $7$ or
$8$ and $\pi\cap \pi(S)=\{2,3\}$. If  $H$ is a Sylow subgroup, then  by \cite[Corollary~C]{ZenkovSylow} it follows that there exist
$x,y\in G$ such that $H\cap H^x\cap H^y=\{e\}$. So assume that  $H$ is a $\{2,3\}$-Hall subgroup of $G$ and $n=5$, $7$
or $8$. In this case  $\Aut(\Alt_n)=\Sym_n$ and, since  $G=HS$, we need to prove the statement of the theorem for the case $G=\Aut(S)=\Sym_n$. Consider all three possibilities for~$n$.

If  $G=\Sym_5$, then $H=\Sym_4$ is a point stabilizer in the natural permutation representation. Since the intersection of every four point stabilizers is trivial and point stabilizers of a transitive group are conjugate, there exist $x,y,z\in G$ such that $H\cap H^x\cap H^y\cap H^z=\{e\}$.

If $G=\Sym_7$, then $H=\Sym_3\times\Sym_4$. Up to conjugation in $G$, sets $\{1,2,3\}$ and $\{4,5,6,7\}$ are $H$-orbits. Direct calculations show that for $x=(2,4)(3,5)$ and $y=(1,2,4)(3,6,5)$ the equality $H\cap H^x\cap H^y=\{e\}$ holds.

If $G=\Sym_8$, then $H=\Sym_4\wr\Sym_2$. Up to conjugation in  $G$ we may assume that $H$ is generated by
$(1,2,3,4),(1,2),(1,5)(2,6)(3,7)(4,8)$. In this case there do not exist elements $x,y,z\in G$ such that $H\cap H^x\cap H^y\cap H^z=\{e\}$. Indeed, the action by right multiplication of  $G$ on the set of right cosets $G:H$ gives an embedding of $G$ into $\Sym_{35}$ (and the image of $G$ under this embedding is primitive). By using \cite{GAP} it is easy to check that on the set  $(G:H)^4$  under this action there exist  $152$ orbits, of order at most $20160=\vert G\vert/2$, i.~e., the stabilizer of every four points is nontrivial. This means that the intersection of every four conjugate with $H$ subgroups (by construction $H$ is a point stabilizer) is nontrivial.

For elements  $x=(4,5)$, $y=(3,5)(4,6)$, $z=(4,7,6,5)$, $t=(2,3,4,5)$ the equality $H\cap H^x\cap H^y\cap H^z\cap H^t=\{e\}$ holds.
We estimate the number of orbits of  $H$ on the set $X$ of four-tuples $$X:=\{(H^x,H^y,H^z,H^t)\mid H\cap H^x\cap H^y\cap H^z\cap
H^t=\{e\}\}.$$ Since the equality $N_G(H)=H$ holds, an element $g\in G$ leaves the four-tuple $(H^x,H^y,H^z,H^t)$ invariant if and only if
$g\in  H^x\cap H^y\cap H^z\cap H^t$. Thus the action of $H$ on $X$ is regular, so the size of every $H$-orbit is equal to $\vert H\vert=1152$. The intersection $H^x\cap H^y\cap H^z\cap H^t$ is a cyclic group of order  $2$, which is generated by $a=(1,8)(2,4)(3,6)(5,7)$. The equality $H\cap
(H^x)^g\cap (H^y)^g\cap (H^z)^g\cap (H^t)^g=\{e\}$ is true if and only if $a^g\not\in H$. More over, the equality  $N_G(H)=H$ implies that four-tuples $$((H^x)^{g_1},(H^y)^{g_1},(H^z)^{g_1},(H^t)^{g_1}),
((H^x)^{g_2},(H^y)^{g_2},(H^z)^{g_2},(H^t)^{g_2})$$ coinside if and only if $g_1g_2^{-1}\in\langle a\rangle$. Direct calculations (or by using  \cite{GAP}) it is easy to check that  $G\setminus H$ contains $72$ conjugate with $a$ elements. Further more  $C_G(a)\simeq
2\wr \Sym_4$ and  $\vert C_G(a)\vert=384$. So the cardinality of  $X$ is equal to $72*192=13824$. Therefore $H$ has at least $13824/1152=12$ orbits on~$X$. 

{\bfseries (II)} Assume that $S$ is either a sporadic group or the Tits group. If $H$ is a Sylow subgroup of $G$ then by
\cite[Corollary~C]{ZenkovSylow} it follows that there exist  $x,y\in G$ such that $H\cap H^x\cap H^y=\{e\}$. So we may assume that the order of  $H$ is divisible by at least two primes.

Assume first that the order of $H$ is divisible by precisely two primes $r,p$ and $H=R: P$, where $R\in \mathrm{Syl}_r(H)$, $P\in\mathrm{Syl}_p(H)$ and $O_p(H)=\{e\}$. Assume also that either  $r$ is odd or  $S=\Aut(S)$. Since for every sporadic group and Tits group the group of its outer automorphisms is a $2$-group, then $R\leq S$ under our assumption. By \cite{MazZen} there exists $x\in S$ such that  $R\cap R^x=\{e\}$. Hence, up to conjugation in $H$ we may assume that $H\cap H^x\leq P$. By \cite{Passman} it follows that there exist $y,z\in R$ such that $P\cap P^y\cap
P^z=\{e\}$. So $(H\cap H^x)\cap (H\cap H^x)^y\cap (H\cap H^x)^z=H\cap H^x\cap H^{xy}\cap H^{xz}=\{e\}$.

If the order of  $H$ is divisible by two primes $r,p$ only and $H=R\times P$, where $R\in \mathrm{Syl}_r(H)$, $P\in\mathrm{Syl}_p(H)$, then one of these primes, say $r$, is odd and again $R\leq S$.  By \cite{MazZen} there exists an element   $x\in S$ such that $R\cap R^x=\{e\}$. More over, by  \cite[Corollary~C]{ZenkovSylow} there exist $y,z\in G$ such that $P\cap P^y\cap P^z=\{e\}$. Hence $H\cap H^x\cap
H^y\cap H^z=\{e\}$.

By \cite[Theorem~6.14 and Table~III]{GrossOdd} and \cite[Theorem~4.1]{RevinDpiOneClass} it follows that a  $\pi$-Hall subgroup $H$ of a sporadic group $S$ such that the order of $H$ is divisible by at least two primes and the structure of $H$ is not covered by cases considered above exists only if either $S\simeq M_{23}$, $H\simeq 2^4:(3\times \Alt_4):2\simeq (2^4:2^2):3^2:2$, or
$S\simeq J_1$, $H\simeq2^3:7:3$. In both cases we have $S=\Aut(S)=G$.

Assume first that   $S\simeq M_{23}$, $H\simeq (2^4:2^2):3^2:2 =R:P:Q$. By \cite{MazZen} there exists $x\in G$ such that
$R\cap R^x=\{e\}$, therefore we may assume that $H\cap H^x\leq P:Q$. By Corollary \ref{IntersectionAbelHallsubgroups} there exists
$y\in R$ such that $(H\cap H^x)\cap (H\cap H^x)^y=H\cap H^x\cap H^{xy}$ and $\vert H\cap H^x\cap H^{xy}\vert\le2$. Again by using Corollary
\ref{IntersectionAbelHallsubgroups}, we find an element $z\in H\cap H^x$ with  $$( H\cap H^x\cap
H^{xy})\cap( H\cap H^x\cap H^{xy})^z=H\cap H^x\cap H^{xy}\cap H^{xyz}=\{e\}.$$

Assume, that $S\simeq J_1\simeq G$, $H\simeq2^3:(7:3)=R: P$, where $R=2^3=O_2(H)\in \mathrm{Syl}_2(H)$ and $P\in \mathrm{Hall}_{\{3,7\}}(H)$, moreover $O_{\{3,7\}}(H)=\{e\}$. By \cite{MazZen} it follows that there exists  $x\in S$ such that  $R\cap R^x=\{e\}$, so, up to conjugation in $H$, we may assume that  $H\cap H^x\leq P$. By \cite{DolfiOdd} there exist $y,z\in R$ such that $P\cap P^y\cap P^z=\{e\}$. Hence $H\cap H^x\cap H^{xy}\cap H^{xz}=\{e\}$.
\qed

Authors

Vdovin Evgeny Petrovitch

Sobolev Institute of mathematics SB RAS

e-mail:vdovin@math.nsc.ru

Zenkov Victor Ivanovitch

Institute of mathematics and mechanics UrB RAS

\end{document}